\documentclass[10pt]{amsart}

\usepackage{graphicx}
\usepackage{marvosym}
\usepackage{ifsym}
\usepackage{lipsum}

\theoremstyle{definition}

\theoremstyle{remark}

\begin{document}

\title[The Gauss-Bonnet formula]{The Gauss-Bonnet formula of a conical metric on a compact Riemann surface}


\author{FANG Hanbing}
\address{Mathematics Department, Stony Brook University, NY 11794, United States}
\email{hanbing.fang@stonybrook.edu}

\author{XU Bin}
\address{CAS Wu Wen-Tsun Key Laboratory of Mathematics and  School of Mathematical \newline \indent Sciences, University of Science and Technology of China, Hefei 230026 China}
\email{bxu@ustc.edu.cn}

\author{YANG Bairui}
\address{School of Mathematical Sciences, University of Science and Technology of China\newline \indent Hefei 230026 China}
\email{bry@mail.ustc.edu.cn}
\date{}

\dedicatory{}

\begin{abstract}
We prove a generalization of the classical Gauss-Bonnet formula for a conical metric on a compact Riemann surface provided that the Gaussian curvature is Lebesgue integrable with respect to the area form of the metric.
We also construct explicitly some conical metrics whose curvature is not integrable.\\

\noindent {\sc Keywords:}  Gauss-Bonnet formula, conical metric, Riemann surface, Gaussian curvature, Lebesgue integrable\\

\noindent {\sc 2000 MR Subject Classification:} 53C21

\end{abstract}

\maketitle
\newcommand\blfootnote[1]{%
\begingroup
\renewcommand\thefootnote{}\footnote{#1}%
\addtocounter{footnote}{-1}%
\endgroup
}


S.-S. Chern \cite{Ch} gave an intrinsic proof of the Gauss-Bonnet formula on a closed orientable Riemannian manfold of even dimension by thinking of the Euler-Poincar\' e characteristic of the manifold as a topological invariant of the tangent sphere bundle to the manifold. M. Atiyah and C. Lebrun \cite[Theorems 2.1-2]{AL} proved modified versions of the Gauss-Bonnet and signature theorems for compact Riemannian 4-manifolds with edge-cone singularities along an embedded 2-manifold.
C. T. McMullen \cite{Mc} presented an analogue of the Gauss-Bonnet theorem for a particular class of stratified spaces, called cone manifolds, by using the notion of spherical complexes.
By using the Green formula, M. Troyanov \cite[Proposition 1]{MR91} generalized the formula to a conical metric on a compact Riemann surface $X$ provided that the Gaussian curvature $K$ of the metric extends to a H\" older continuous function on $X$. Some years ago, Q. Chen asked the second author whether the H\" older continuity of $K$ could be replaced by the $L^1$ integrability of $K$ with respect to the area form
of the metric.  We shall give an affirmative answer to his question.

Before stating our theorem in a precise form, we at first introduce the notations relevant to conical metrics on Riemann surfaces. Let $p_1,\cdots, p_n$ be $n$ distinct points on a compact Riemann surface $S$ with genus $g_S\geq 0$ and $\beta_1,\cdots,\beta_n$
be real numbers lying in $(-1,\,\infty)$. We say that ${\rm d}s^2$ is a {\it conical metric representing the real divisor $D=\sum_{j=1}^n\, \beta_j\,p_j$} on $S$ if and only if it satisfies the following two properties:

\blfootnote{\noindent{\bf \hspace{-0.55cm} Foundation Item:} Supported in part by the National Natural Science Foundation of China (Grant Nos. 12271495, 11971450 and 12071449) and the CAS Project for Young Scientists in Basic Research (YSBR-001).\\
\noindent{\bf Biographies:} FANG Hanbing(1998.10-), male, native of Chizhou, Anhui,  PhD student of Stony Brook University, engages in Differential Geometry; XU Bin(1976.08-), male, native of Yiyang,
Hunan, associate professor of University of Science and Technology of China, engages in Differential
Geometry; YANG Bairui(1998.05-, corresponding author), native of Nanjing, Jiangsu, master student of University of Science and Technology of China, engages in Geometric Analysis.}

\begin{enumerate}
\item[(i)] ${\rm d}s^2$ is a smooth conformal metric on $S\backslash {\rm supp}\, D=S\backslash\{p_1,\cdots, p_n\}$, whose area form we denote by ${\rm d}A$;
\item[(ii)] for each $1\leq j\leq n$, there exists a complex coordinate chart $(U_j,\, z_j)$ centered at $p_j$ such that $z_j(p_j)=0$, ${\rm d}s^2=e^{2\varphi_j(z_j,\,\overline{z_j})}\,|{\rm d}z_j|^2$ on $U_j\backslash \{p_j\}$
    and the smooth function $\varphi_j-\beta_j\,\ln\,|z_j|$ on $U_j\backslash \{z_j=0\}$ extends {\it continuously} to $z_j=0$. In particular, the area form ${\rm d}A$ of ${\rm d}s^2$ equals $\frac{\sqrt{-1}}{2}e^{2\varphi_j}{\rm d}z_j\wedge {\rm d}\overline{z_j}$ on $U_j\backslash \{z_j=0\}$. Note that the area of ${\rm d}s^2$ on $S$ is finite.
\end{enumerate}
At this time, ${\rm d}s^2$ has conical angle $2\pi(1+\beta_j)$ at $p_j$.
We call $\deg\,D=\sum_{j=1}^n\,\beta_j$ the {\it degree} of $D$ and $\chi(S,\,D)=2-2g_S+\deg\,D$
the {\it singular Euler number} of pair $(S,\,D)$. We denote by $K$ the Gaussian curvature of ${\rm d}s^2$ on
the punctured surface $S\backslash {\rm supp}\, D$.
Troyanov's generalization of the Gauss-Bonnet formula
says that
\begin{equation}
\label{equ:GB}
\frac{1}{2\pi}\, \int_{S\backslash {\rm supp}\, D}\, K\,{\rm d}A=\chi(S,\, D).
\end{equation}
provided that $K$ extends to a H\" older continuous function on $S$ (\cite[Proposition 1]{MR91}).\\

\noindent{\bf Theorem} {\it The Gauss-Bonnet formula \eqref{equ:GB} holds for the conical metric ${\rm d}s^2$ representing $D$ on $S$ if and only if
its Gaussian curvature $K$ is $L^1$ integrable with respect to the area form ${\rm d}A$ on $S\backslash {\rm supp}\, D$.}\\

We need the following lemma to prove the theorem.\\


\noindent{\bf Lemma} {\it
Let $\mathrm{d}s^{2}=e^{2u}\vert z\arrowvert ^{2\beta}\arrowvert \mathrm{d}z\vert ^{2}$ be a conical metric on the disk $U= \{|z|<1 \}$ which has conical angle $2\pi(1+\beta)>0$ at $z=0$ and whose Gaussian curvature $K$ is $L^1$ integrable with respect to the area form ${\rm d}A=\frac{\sqrt{-1}}{2}\,e^{2u}\vert z\arrowvert ^{2\beta}\,{\rm d}z\wedge {\rm d}\bar z$ on $U^*=\{0<|z|<1\}$. Then we have
$$\liminf_{\epsilon\to 0^+}\,\int_{|z|=\epsilon}\left(\left| \frac{\partial u}{\partial z}\right|+\left|\frac{\partial u}{\partial \bar{z}}\right|\right)|\mathrm{d}z|=0.$$
}
\begin{proof} We prove the lemma in three steps.
\begin{enumerate}

\item[(i)] We show at first that {\it the equality $\displaystyle{\Delta u:=4\frac{\partial ^{2}u}{\partial z \partial \bar{z}}=-K\vert z\vert ^{2\beta}e^{2u}}$ holds in the sense of distribution in $U$}. Actually, by the very definition of Gaussian curvature, the equality holds point-wisely on $U\backslash \{z=0\}$.  Moreover, by the definition of conical metric, $u$ is continuous on $U$ and $\Delta u$ is a distribution on $U$. Since $K$ is $L^1$ integrable with respect to ${\rm d}A$, so is $h:=-K\vert z\vert^{2\beta}e^{2u}$ with respect to the Euclidean area form ${\rm d}\bar{A}=\frac{\sqrt{-1}}{2}{\rm d}z\wedge{\rm d}\bar z$, and $h$ is also
    a distribution in $U$. Hence, the distribution
  $\Delta u-h$ is supported at $z=0$ and equals some finite sum $\sum C_{\alpha}\partial^{\alpha}\delta_{0}$ by \cite[pp.46-47, Theorem 2.3.4]{LH90}. We take an arbitrary multi-index $\alpha$ and choose $\phi \in C_{0}^{\infty}(U)$ such that $\partial^{\alpha}\phi (0)\neq 0$ and $\partial^{\beta}\phi (0)= 0$ for all $\beta \neq \alpha$. Denote $\phi_{k}(\cdot)=\phi(k\cdot)$. Since
$h$ is an $L^1$ function with respect to $\mathrm{d} \bar{A}$, we have
$$ |h(\phi_{k})| = \left| \int_{U}h\phi_{k} \,\mathrm{d}\bar{A}\right| \leq \sup |\phi_{k}|\int_{{\rm supp}\,\phi_{k}}|h|\,\mathrm{d}\bar{A} \to 0 \quad {\rm as}\quad k\to \infty.$$ \\
Since $u$ is continuous at $z=0$ and $\displaystyle{\int_{U}\Delta \phi = 0}$, $|\Delta u(\phi_{k})|$ equals
$$\left|\int_{U}u\left(\frac{\cdot}{k}\right) \Delta \phi\, \mathrm{d}\bar{A} - \int_{U}u(0)\Delta \phi \,\mathrm{d}\bar{A} \right| \leq \int_{U}\left|u\left(\frac{\cdot}{k}\right)-u(0)\right||\Delta \phi| \,\mathrm{d}\bar{A} \to 0\,\,{\rm as}\,\,k\to\infty.
$$
Therefore we obtain $(\Delta u -h)(\phi_{k})$ tends to $0$ as $k$ goes to $\infty$. On the other hand, $(\Delta u -h)(\phi_{k})  = C_{\alpha} k^{|\alpha|} \partial^{\alpha}\phi(0)$, which implies $C_{\alpha} = 0$. Thus $\Delta u = h$ as distributions in $U$. \\

\item[(ii)] We show that {\it the estimate $\displaystyle{\int_{0<|z|<r}|\nabla u|^{2}\,\mathrm{d}\bar{A}<\infty}$ holds as $0<r<1$, where $\displaystyle{|\nabla u|^{2}=4\left(\left| \partial u/\partial z\right|^2+\left|\partial u/\partial \bar{z}\right|^2\right)}$}. To this end, we only prove the case of $r=1/2$ for this estimate for simplicity. At first we choose $\{\chi_{\epsilon}:0<\epsilon<1\}\subset C_0^\infty({\Bbb C})$  such that $\chi_\epsilon\geq 0$, $\int_{\mathbb{C}}\chi_{\epsilon}=1$ and ${\rm supp}\, \chi_{\epsilon}\subset \{|z|\leq \epsilon\}$.
We could further assume ${\rm supp}\, u\subset \{|z|<3/4\}$ by a suitable cut-off function if necessary.
Since $u$ is continuous on ${\Bbb C}$, the sequence of convolutions $\{u_{k}:=\chi_{\frac{1}{k}}*u\}\subset C_0^\infty({\Bbb C})$ converges uniformly to $u$ on ${\Bbb C}$ as $k$ goes to $\infty$. Moreover, since $u$ is smooth in ${\Bbb C}\backslash\{0\}$, $|\nabla u_{k}|^{2}$ converges to $|\nabla u|^{2}$ uniformly in any compact subset of ${\Bbb C}\backslash\{0\}$ as $k$ goes to $\infty$. By the Fatou lemma, we have
$$\int_{0<|z|<1/2}|\nabla u|^{2}  \mathrm{d}\bar{A}\leq \liminf_{k\to \infty}\int_{|z|<1/2}|\nabla u_{k}|^{2} \mathrm{d}\bar{A}.$$
Hence the problem is reduced to showing that {\it $\displaystyle{\int_{|z|<1/2}|\nabla u_{k}|^{2}\mathrm{d}\bar{A}}$ is uniformly bounded from above}. By using integration by part, we have
$$\int_{|z|<1/2}|\nabla u_{k}|^{2} \mathrm{d}\bar{A}=-\int_{|z|<1/2} u_{k}\Delta u_{k}\mathrm{d}\bar{A}+\int_{|z|=1/2} u_{k}\frac{\partial u_{k}}{\partial \vec {n}}\,|{\rm d}z|.$$
Since $\displaystyle{u_k\frac{\partial u_{k}}{\partial \vec {n}}}$ converges to 
$\displaystyle{u\frac{\partial u}{\partial \vec {n}}}$ uniformly on $|z|=1/2$ as $k$ goes to $\infty$, the second summand on the right hand side is uniformly bounded. For the first one,
as $k\geq 4$, we have
\begin{align*}
  \int_{|z|<1/2}|\Delta u_{k}|\mathrm{d}\bar{A}&=\int_{|z|<1/2}|\Delta (\chi_{1/k}\ast u)|\mathrm{d}\bar{A} =\int_{|z|<1/2}|\chi_{1/k}\ast(\Delta u)|\mathrm{d}\bar{A}\\
  &=\int_{|z|<1/2}|\chi_{1/k}\ast h|\mathrm{d}\bar{A}=\int_{|z|<1/2}\left|\int_{|\hat{z}|
  <\frac{1}{2}+\frac{1}{k}}\chi_{\frac{1}{k}}(z-\hat{z})h(\hat z)\right|\\
  &\leq \int_{\mathbb{C}}\chi_{\frac{1}{k}}\int_{|\hat z|<\frac{1}{2}+\frac{1}{k}} |h|
  \leq\int_{|\hat z|<3/4}\,|h|<\infty,
\end{align*}
where we use the statement in Step (i) in the third equality. Hence,
we are done since $u_k$ is uniformly bounded on ${\Bbb C}$.

\item[(iii)] At last, we prove the lemma by contradiction. Suppose that there
exists some positive constant $C$ such that when $0<\epsilon<1$ is sufficiently small,
$$\int_{|z|=\rho }\left(\left|\frac{\partial u}{\partial z}\right|+\left|\frac{\partial u}{\partial \bar z}\right|\right)\mathrm{d}\theta\geq \frac{C}{\rho}\quad \text{for all}\quad 0<\rho<\epsilon.$$
Then we have
\begin{eqnarray*}
\int_{0<|z|<\epsilon } 2\left(\left|\frac{\partial u}{\partial z}\right|^{2}+\left|\frac{\partial u}{\partial \bar z}\right|^{2}\right)\mathrm{d}\bar{A}
&\geq&\int_{0}^{\epsilon} \rho\,{\rm d}\rho \int_{0}^{2\pi}\,\left(\left|\frac{\partial u}{\partial z}\right|+\left|\frac{\partial u}{\partial \bar z}\right|\right)^2 \mathrm{d}\theta\\
 ({\rm Cauchy-Schwarz})\quad &\geq&  \int_{0}^{\epsilon}\,\rho\,{\rm d}\rho\, \frac{1}{2\pi}\cdot\left(\int_{|z|=\rho}\left(\left|\frac{\partial u}{\partial z}\right|+\left|\frac{\partial u}{\partial \bar z}\right|\right)\mathrm{d}\theta\right)^2 \\
&\geq&\int_{0}^{\epsilon}\,\frac{\rho}{2\pi}\,\left(\frac{C}{\rho}\right)^2\,{\rm d}\rho=\infty,
\end{eqnarray*}
which contradicts the estimate proved in the second step.

\end{enumerate}

\end{proof}

The proof of the theorem in what follows goes similarly as that of \cite[Proposition 1]{MR91} with  necessary modifications.

\begin{proof} We only need to show the sufficient part since the necessary one is clear.
 Choose a smooth conformal metric $\mathrm{d}s_{1}^{2}$ which has Gaussian curvature  $K_1$ and area form $\mathrm{d}A_1$ on $S$. By the classical Gauss-Bonnet formula, we have
$\displaystyle{\frac{1}{2\pi}\int_{S} K_{1} \mathrm{d}A_{1}=\chi (S)}$.
Moreover, there exists a smooth real function $v$ on $S\backslash {\rm supp}\, D$ such that $\mathrm{d}s^{2}=e^{2v}\mathrm{d}s_{1}^{2}$ and
$\displaystyle{K\mathrm{d}A=K_{1}{\rm d}A_1-\mathrm{d}*\mathrm{d}v}$,
where $\ast$ is the conjugation operator acting one-forms on $S$ , and locally it has form $\displaystyle{*\mathrm{d}v=-i\frac{\partial v}{\partial z}\mathrm{d}z +i\frac{\partial v}{\partial \bar {z}}\mathrm{d}\bar{z}}$ (\cite[p. 25, (3.7.1)]{FK}).
It suffices to show the equality
$\displaystyle{-\frac{1}{2\pi}\int _{S\backslash {\rm supp}\, D}\mathrm{d}*\mathrm{d}v=\deg\, D}$.
Let $z_j$ be complex coordinate centered at $p_j$ and $D_{j}(\delta)=\{|z_j|<\delta\}$  mutually disjoint open disks of radius $\delta$ centered at $p_j$ for all $1\leq j\leq n$.
Recall that $p_j$ is a conical singularity with angle $2\pi(1+\beta _{j})$, $K\mathrm{d}A=K_{1}{\rm d}A_1-\mathrm{d}*\mathrm{d}v$ and $K$ is $L^1$ integrable with respect to ${\rm d}A$ on $S\backslash{\rm supp}\, D$. Then  we have $v=\beta _{j} \log |z_j|+u_j$, where $u$ is smooth on $\{0<|z_j|<\delta\}$ and $2\frac{\partial^2 u_j}{\partial z_j \partial \overline{z_j} }$ is $L^1$ integrable with respect to $\frac{\sqrt{-1}}{2}{\rm d}z_j\wedge {\rm d}\overline{z_j}$ on $0<|z_j|<\delta$.
Since
\begin{equation*}
\frac{1}{2\pi}\int _{\partial D_{j}(\delta)}*\mathrm{d}v=\frac{\beta _{j}}{2\pi}\int _{\partial D_{j}(\delta)} *\mathrm{d}\log |z_j|+\frac{1}{2\pi}\int _{\partial D_{j}(\delta)} *\mathrm{d}u_j=\beta _{j}+\frac{1}{2\pi}\int _{\partial D_{j}(\delta)}*\mathrm{d}u_j,
\end{equation*}
by the Green formula, we have
\begin{eqnarray*}
-\frac{1}{2\pi}\int _{S \setminus \bigcup_{j=1}^n\, D_{j}(\delta)} \mathrm{d}*\mathrm{d}v=\frac{1}{2\pi} \sum _{j=1}^{n}\int_{\partial D_{j}(\delta)}*\mathrm{d}v
=\deg\, D+\sum_{j=1}^n\,\frac{1}{2\pi}\int _{\partial D_{j}(\delta)}*\mathrm{d}u_j.
\end{eqnarray*}
We complete the proof of the theorem by taking the lower limit of both the left and right hand sides as $\delta\to 0^+$ and using the lemma.

\end{proof}

At last we give some explicit conical metrics whose Gaussian curvature is {\it not} $L^1$ integrable and whose angles take all positive numbers  in the following: \\

\noindent
{\bf Example} {\it Since the function $\displaystyle{\phi(r):= \alpha\left(r+r^3\,\sin\frac{1}{r}\right)}$ is always positive in the interval $(0,\,1)$ for each $\alpha>0$, $\displaystyle{\mathrm{d}s^{2} = \mathrm{d}r^{2}+\phi(r)^{2} \mathrm{d}\theta^2}$ is a rotational symmetry metric under the polar coordinate chart $(r,\,\theta)\in (0,\,1)\times [0,\, 2\pi]$, of which the point of $r=0$ might be a singularity.
Then the Gaussian curvature $K$ of ${\rm d}s^2$
is not $L^1$ integrable with respect to the area form $\phi(r){\rm d}r\wedge {\rm d}\theta$ near point $r=0$.
Moreover, under some complex coordinate $z=z(r,\,\theta)$, ${\rm d}s^2$ is a smooth conformal metric
in $0<|z|<1$ such that  $z=0$, i.e. $r=0$, is a conical singularity of angle $2\pi\alpha$.}

\begin{proof}
Since the Gauss curvature $\displaystyle{K = - \frac{{\rm d}^2\phi/{\rm d}r^2}{\phi(r)}}$,  for each $0<\rho<1$, we have
\begin{eqnarray*}
\int_0^\rho\,|K|\phi(r)\mathrm{d}r &=& \int_0^\rho\, \left|\frac{{\rm d}^2\phi}{{\rm d}r^2}\right|\,\mathrm{d}r=
\alpha\,\int_0^\rho\, \left|6r\sin\frac{1}{r} - 4\cos\frac{1}{r} - \frac{1}{r}\sin\frac{1}{r}\right|\, {\rm d}r\\
&\geq& \alpha\,\int_0^\rho\, \left|\frac{1}{r}\sin\frac{1}{r}\right|\, {\rm d}r-
\alpha\,\int_0^\rho\, \left(\left|6r\sin\frac{1}{r}\right|+4\left|\cos\frac{1}{r}\right|\right)\, {\rm d}r\\
&=&\infty\,.
\end{eqnarray*}
This shows that $K$ is not $L^1$ integrable in any neighborhood of $r=0$. Defining
\[z:=e^w\quad {\rm where}\quad  w:=-v+\sqrt{-1}\theta \quad {\rm and}\quad v=v(r):=\int_r^1\,\frac{{\rm d}\rho}{\phi(\rho)},\]
we find that $z$ could take all non-zero complex values of modulus less than $1$ as $(r,\,\theta)\in (0,\,1)\times [0,\, 2\pi]$.
Since $\displaystyle{\lim_{r\to 0^+}\, \frac{\alpha v(r)}{\ln\,(1/r)}=1}$,  we have as $z\to 0$
\begin{eqnarray*}
{\rm d}s^2&=&\phi^2\big(\mathrm{d}v^{2}+\mathrm{d}\theta^2\big)=\phi^2|{\rm d}w|^2=\frac{\phi^2}{|z|^2}|{\rm d}z|^2\sim \alpha^2|z|^{2\alpha-2}|{\rm d}z|^2,
\end{eqnarray*}
i.e. $z=0$ is a conical singularity of angle $2\pi\alpha$.

\end{proof}

\noindent\textbf{Acknowledgement:}
B.X. would like to express his deep gratitude to Professor Qing Chen at USTC for his stimulating question.
The authors thank the anonymous referee for her/his constructive comments which improve the manuscript greatly.

\bibliographystyle{amsplain}

\end{document}